\documentclass{amsart}
\usepackage{amsfonts,amssymb,amsmath}
\usepackage{url}

\usepackage{graphicx} 

\newtheorem{thrm}{Theorem}[section]

\newtheorem{cor}[thrm]{Corollary}
\theoremstyle{definition}

\newtheorem{example}[thrm]{Example} 
\theoremstyle{remark} 
\newtheorem*{remark}{Remark} 

\newcommand{\C}{\mathbb{C}}
\newcommand{\D}{\mathbb{D}}
\newcommand{\E}{\mathbb{E}}
\newcommand{\Half}{\mathbb{H}}
\newcommand{\Halft}{\mathbb{H}_{\,t}}
\newcommand{\R}{\mathbb{R}}
\newcommand{\Z}{\mathbb{Z}}
\newcommand{\A}{\mathcal{A}}
\newcommand{\FominEvent}{\mathcal{C}}
\newcommand{\curves}{\mathcal{K}}
\newcommand{\LE}{\mathcal{L}}
\newcommand{\n}{\mathbf{n}}
\newcommand{\Prob}{\mathbf{P}}
\newcommand{\x}{\mathbf{x}}
\newcommand{\y}{\mathbf{y}}
\newcommand{\w}{\omega}
\newcommand{\eps}{\varepsilon}
\newcommand{\bd}{\partial}  
\newcommand{\ds}{\displaystyle}
\newcommand{\SLE}{\operatorname{SLE}}
\newcommand{\hcap}{\operatorname{hcap}} 
\newcommand{\dist}{\operatorname{dist}}

\begin{document}

\title[The scaling limit of Fomin's identity]{The scaling limit of Fomin's identity for two paths in the plane}

\author{Michael J.~Kozdron}
\address{
University of Regina\\
Department of Mathematics \& Statistics\\
College West 307.14\\
Regina, SK S4S 0A2}
\email{kozdron@stat.math.uregina.ca}
\thanks{Research supported by the Natural Sciences and Engineering Research Council (NSERC) of Canada.}

\begin{abstract}
We review some recently completed research that establishes the scaling limit of Fomin's identity for loop-erased random walk on $\Z^2$ in terms of the chordal Schramm-Loewner evolution (SLE) with parameter $2$. In the case of two paths,
we provide a simplified proof of the identity for loop-erased random walk and simple random walk, and prove directly that the corresponding identity holds for chordal $\SLE_2$ and Brownian motion.  We also include a brief introduction to SLE and discussion of the relationship between $\SLE_2$ and loop-erased random walk.\\

\noindent\textsc{Resum\'e}. Nous passons en revue de la recherche r\'ecemment r\'ealis\'ee qui \'etablit la limite de l'identit\'e de Fomin pour la marche al\'eatoire \`a boucles effac\'ees sur $\Z^2$ en termes du processus Schramm-Loewner (o\`u SLE pour Schramm-Loewner evolution) avec param\`etre 2.
Dans le cas de deux chemins, nous fournissons une preuve simplifi\'ee de l'identit\'e pour la marche al\'eatoire \`a boucles effac\'ees et la marche al\'eatoire simple, et prouvons ordonner que l'identit\'e correspondante se tient pour $\SLE_2$ et le mouvement brownien.
Nous incluons \'egalement une br\`eve introduction au processus Schramm-Loewner et une discussion du rapport entre $\SLE_2$ et la marche al\'eatoire \`a boucles effac\'ees.
\end{abstract}

\maketitle


\section{Introduction}

The primary purpose of this paper is to review some recently completed research that has established the scaling limit of Fomin's identity for loop-erased random walk on $\Z^2$ in terms of the chordal Schramm-Loewner evolution (SLE) with parameter $\kappa=2$. For the complete details, including extensions of these results, consult the original papers~\cite{Fom1},~\cite{Koz},~\cite{KozL} and~\cite{KozL2}. We have decided to discuss the case of $n=2$ paths exclusively. This choice is partly pedagogical, and it is our hope that the reader will find the particular special cases discussed in the present work to be useful in understanding the general results of the original papers. Furthermore, explicit calculations can be performed in the case of $n=2$ paths, and this choice allows us to present simplified proofs of these theorems in this case. It must also be noted that
we shall only discuss two-dimensional results at present. 
Therefore, we will consider $\R^2 \cong \C$, and will write any of $w$, $x$, $y$, $z$ to denote points in $\C$. A simple random walk on $\Z^2$ will be denoted by $S_j$, $j=0, 1, 2, \ldots$, and $B_t$, $t \ge 0$,  will denote a complex Brownian motion. When a one-dimensional Brownian motion is needed, we will write it as $\{W_t, t\ge 0\}$.

It is assumed that the reader has an understanding of random walk and Brownian motion. Although some familiarity with SLE would be helpful,  it is not necessary, and in order for this paper to be as widely accessible as possible, Section~\ref{SLEintrosec} provides a brief introduction to SLE  and discusses the relationship between $\SLE_2$ and loop-erased random walk. The outline of the remainder of the paper is as follows.  Section~\ref{motivationsec} provides some motivation for this paper. Section~\ref{LERWsec} discusses only discrete results including a review of the definition of loop-erased random walk and Fomin's identity.  Section~\ref{epksect} reviews the excursion Poisson kernel, and finally in Section~\ref{nonintsec} we compute the non-intersection probability of $\SLE_2$ and Brownian motion, and explain how it is the natural continuous analogue of Fomin's identity. 


\section{Motivation}\label{motivationsec}

This paper (and~\cite{KozL2} more generally) is the result of the answer to the following question. Suppose that $x$ and $y$ are real numbers with $0<x<y<\infty$.  What is the probability that a chordal $\SLE_2$ from $0$ to $\infty$ in the upper half plane $\Half$ and a Brownian motion excursion from $x$ to $y$ in $\Half$ do not intersect? (See Section~\ref{SLEintrosec} for a brief introduction to SLE and a discussion of the relationship between $\SLE_2$ and loop-erased random walk.) The motivation for asking this question is that the probability under consideration is the natural continuous analogue of the probability that arises in Fomin's identity. (See Section~\ref{LERWsec} for a review of loop-erased random walk and Fomin's identity.) In fact, Fomin's original identity~\cite{Fom1} expressed a particular ``crossing probability'' for loop-erased random walk in terms of the determinant of the hitting matrix for simple random walk, and in that work he conjectured that this identity holds for continuous processes:
\begin{quote}\label{quote}
``\dots we do not need the notion of loop-erased Brownian motion. 
Instead, we discretize the model, compute the probability, and then pass to the limit."
\end{quote}
It is well-known that Brownian motion is the scaling limit of random walk. In~\cite{KozL}, the technical details necessary to complete this conjectured program were first carried out, and the scaling limit of the determinant of the hitting matrix for simple random walk was shown to be the  determinant of the hitting matrix for Brownian motion.
Since the scaling limit of loop-erased random walk is known~\cite{LSW} to be $\SLE_2$, it is natural to ask if the scaling limit of the ``crossing probability'' for loop-erased random walk can be given directly in terms of an $\SLE_2$ probability, and if that $\SLE_2$ probability is equal to the determinant of the hitting matrix for Brownian motion. The answer is in the affirmative as Theorem~\ref{FI_SLE} shows explicitly in the case of two paths.

Indeed no notion of loop-erased Brownian motion was needed!  It should be noted  that although the notion of loop-erased Brownian motion is not well-defined, there is a sense in which $\SLE_2$ can be thought of as Brownian motion without loops.  This description is given in terms of the Brownian loop soup: adding Brownian loops to an $\SLE_2$ path is one way to produce a Brownian motion.  This result is not relevant for the present paper, but the interested reader can consult~\cite{rwloopsoup} and~\cite{loopsoup} for more precise statements.


\section{Fomin's identity for loop-erased random walk}\label{LERWsec}

Suppose that $A \subset \Z^2$. We define the \emph{(outer) boundary} of $A$ to be $\bd A := \{z \in \Z^2 \setminus A : \dist(z,A)=1\}$, and we say that $A$ is \emph{simply connected} if both $A$ and $\Z^2 \setminus A$ are non-empty and connected. 
Let $\A$ denote the collection of simply connected subsets $A$ of $\Z^2$ containing the origin. Let $S_j$, $j=0,1,\ldots$, denote two-dimensional simple random walk, and suppose that $\tau_A := \inf \{j \ge 1: S_j \not\in A\}$.
We say that a path $\w := [\w_0, \ldots, \w_k]$ is a \emph{discrete excursion in $A$} if  $\w_0$, $\w_k \in \bd A$; $\w_1, \ldots, \w_{k-1} \in A$; and $|\w_j-\w_{j-1}|=1$, $j=1, \ldots, k$. The \emph{length} of $\w$ is $|\w|=k$; it is implicit that $2 \le k <\infty$. Finally, we write $\curves_A$ for the set of discrete excursions in $A$, and define the \emph{simple random walk excursion measure on $A$} to be the measure on $\curves_A$ which gives mass $4^{-k}$ to each discrete excursion in $A$ of length $k$. Note that the excursion measure of $\w$ is the probability that the first $k$ steps of a simple random walk starting at $\w_0$ are the same as $\w$. If $\w$ is a discrete excursion in $A$, let 
$p(\w) := \Prob^{\w_0}\{S_j=\w_j, \; j=0, \ldots, |\w| \}$.
 If $z \in A$, $y \in \bd A$, let the \emph{discrete Poisson kernel} 
 $h_A(z,y)$  be the probability that a simple random walk starting
at $z$ leaves $A$ at $y$; that is, $h_A(z,y):=\Prob^z\{S_{\tau_A}=y\}$.
If $x$, $y \in \bd A$,  let the \emph{discrete excursion Poisson kernel} $h_{\bd A}(x,y)$ be the probability
that a simple random walk starting at $x$ takes its first step
into $A$ and then leaves $A$ at $y$.
That is, 
\begin{equation}\label{feb22eq1}
h_{\bd A}(x,y) := \Prob^x\{S_{\tau_A}=y, \; S_1\in A\} = \sum_{\w \in \curves_A(x,y)} p(\w)
\end{equation}
where we  write $\curves_A(x,y)$ to denote the set of discrete excursions in $A$ with endpoints $x$, $y \in \bd A$.

We now briefly review 
the definition of the loop-erased random walk; see~\cite[Chapter 7]{LawlerGreen}
 for more
details. Since simple random walk on $\Z^2$ 
is recurrent, it is not possible to construct loop-erased random walk by 
erasing loops from an infinite walk.  However, the following loop-erasing 
procedure makes perfect sense since it assigns to each finite simple random walk
 path a self-avoiding path.  Let $S := [S_0, S_1, \ldots, S_k]$ be a
simple random walk path of length $k$. We construct $\LE(S)$, the loop-erased
 part of $S$, recursively as follows. If $S$ is already self-avoiding, 
set $\LE(S):=S$.  Otherwise, let $s_0 := \max\{j : S_j=S_0\}$, and for $i > 0$,
 let $s_i := \max\{j : S_j = S_{s_{i-1}+1} \}$. If we let $m := \min\{i : s_i=k\}$, then 
$\LE(S) := [S_{s_0}, S_{s_1}, \ldots, S_{s_m}]$.

 Suppose that $A \in \A$ and $x^1,\ldots,x^n,y^n,
\ldots,y^1$ are distinct points in $\bd A$, ordered 
counterclockwise. For $i=1, \ldots, n$, 
let $\LE^i := \LE(S^i)$ be the loop erasure of the 
path $[S^i_0 = x^i, S^i_1, \ldots, S^i_{\tau^i_A}]$, 
and let $\FominEvent := \FominEvent(x^1, \ldots, x^n, y^n, 
\ldots, y^1; A)$ be the event that both
\begin{equation}  \label{dec14.1}
 S^i_{\tau^i_A} = y^i, \quad i=1, \ldots, n,
\end{equation}
 and
\begin{equation}  \label{dec14.2}
 S^i[0, \tau^i_A]\cap (\LE^1 \cup \cdots \cup \LE^{i-1})=
\emptyset, \quad i=2,\ldots, n.
\end{equation}

In 2001, S.~Fomin~\cite{Fom1} proved the following identity which relates the determinant of a matrix of simple random walk probabilities to a ``crossing probability'' for loop-erased random walk. 
 
\begin{thrm}[Fomin's Identity]\label{fomintheorem}
If $\FominEvent$ is the event defined above, and 
\begin{equation*}
\mathbf{h}_{\bd A}(\mathbf{x},\mathbf{y}) := 
\begin{bmatrix}
h_{\bd A}(x^1, y^1) &\cdots & h_{\bd A}(x^1, y^n) \\
\vdots      &\ddots &\vdots       \\
h_{\bd A}(x^n, y^1) &\cdots & h_{\bd A}(x^n, y^n)
\end{bmatrix}
\end{equation*}
where $\mathbf{x} := (x^1, \ldots, x^n)$, $\mathbf{y} := (y^1, \ldots, y^n)$, 
then $\Prob(\FominEvent) = \det \mathbf{h}_{\bd A}  (\mathbf{x},\mathbf{y})$.
\end{thrm}

\begin{remark}
We note that the conditional probability that~(\ref{dec14.2})
holds given~(\ref{dec14.1}) holds is
\begin{equation}\label{conditionalFomin}
   \det \left[\frac{h_{\bd A}(x^i,y^{\ell})}{h_{\bd A}(x^i,y^i)}\right]_{1 \leq 
   i,\ell \leq n} 
=\frac{ \ds \det\mathbf{h}_{\bd A}  (\mathbf{x},\mathbf{y}) }{\ds \prod_{i=1}^n h_{\bd A}(x^i,y^i)}. 
\end{equation}
The first approach taken to derive a scaling limit of Fomin's identity and establish the conjecture given in Section~\ref{motivationsec} was to show that~(\ref{conditionalFomin}) converged to the appropriate Brownian motion quantity as the lattice spacing $\delta \to 0$.
This was first accomplished in~\cite{KozL}, and is briefly discussed at the end of Section~\ref{epksect}.
\end{remark}

We end this section with the specific case of two paths for which a simpler proof can be given by ``counting sample paths.'' 

\begin{thrm}[Fomin's Identity for LERW (version for two paths)]\label{FI_LERW_b}
Suppose that $A \in \A$ and $x^1$, $x^2$, $y^2$, $y^1$ are four points ordered counterclockwise around $\bd A$.  If $\LE^1$ is the path of a loop-erased random walk excursion from $x^1$ to $y^1$, and
 $S^2$ is the path of a simple random walk excursion from $x^2$ to $y^2$, then
\begin{equation}\label{fomineq1}
\Prob\{\,
\LE^1 \cap S^2 =\emptyset \,\} =\frac{ \ds \det\mathbf{h}_{\bd A}  (\mathbf{x},\mathbf{y}) }{h_{\bd A}(x^1,y^1)\, h_{\bd A}(x^2,y^2)}.
\end{equation}
\end{thrm}

\begin{proof}
As noted earlier,~(\ref{fomineq1}) represents the conditional probability that the loop-erasure of a first simple random walk excursion starting from $x^1$ and a second simple random walk excursion starting from $x^2$ do not intersect given that the first simple random walk exits at $y^1$ and the second simple random walk exits at $y^2$. Therefore, the key step in proving this theorem is to show that if
\begin{equation}\label{fomeq}
q := q(x^1,x^2,y^2,y^1; A) := \sum p(\w^1)p(\w^2)
\end{equation}
where the sum is over all $\w^1 \in \curves_A(x^1,y^1)$ and $\w^2 \in \curves_A(x^2,y^2)$ with $\LE(\w^1) \cap \w^2 = \emptyset$, then
$q = \det \mathbf{h}_{\bd A}(\mathbf{x},\mathbf{y})
$.
We know from~(\ref{feb22eq1}) that
$$h_{\bd A}(x^i,y^j) := \sum_{\w^{ij} \in \curves_A(x^i,y^j)} p(\w^{ij}),$$
and so we have
\begin{equation}\label{fomeq2}
\det
\begin{bmatrix}
h_{\bd A}(x^1,y^1) &h_{\bd A}(x^2,y^1) \\
h_{\bd A}(x^1,y^2) &h_{\bd A}(x^2,y^2) \\
\end{bmatrix}
=
\sum_{\w^{11}} \sum_{\w^{22}} p(\w^{11}) p(\w^{22})
-
\sum_{\w^{12}} \sum_{\w^{21}} p(\w^{12}) p(\w^{21}).
\end{equation}
Let $\Gamma^1$ denote the set of ordered pairs $(\w^{11},\w^{22})
 \in \curves_A(x^1,y^1) \times \curves_A(x^2,y^2)$ such that $\LE(\w^{11}) \cap \w^{22} \neq \emptyset$ so that
 $$ q = \sum_{(\w^{11},\w^{22}) \in  \curves_A(x^1,y^1) \times \curves_A(x^2,y^2) \setminus \Gamma^1} p(\w^{11}) p(\w^{22}).$$
Let $\Gamma^2$ denote the set of ordered pairs $(\w^{12},\w^{21})\in \curves_A(x^1,y^2) \times \curves_A(x^2,y^1)$, and note that 
$\LE(\w^{12}) \cap \w^{21} \neq \emptyset$ for every $(\w^{12},\w^{21})\in \curves_A(x^1,y^2) \times \curves_A(x^2,y^1)$.  Thus, we can express~(\ref{fomeq2}) as
$$
\det
\begin{bmatrix}
h_{\bd A}(x^1,y^1) &h_{\bd A}(x^2,y^1) \\
h_{\bd A}(x^1,y^2) &h_{\bd A}(x^2,y^2) \\
\end{bmatrix}
= 
q+
\sum_{\Gamma^1} p(\w^{11}) p(\w^{22})
-
\sum_{\Gamma^2} p(\w^{12}) p(\w^{21})
.$$
We will now show that there exists a one-to-one correspondence between $\Gamma^1$ and $\Gamma^2$, denoted by $(\w^{11},\w^{22}) \leftrightarrow (\Lambda\w^{11}, \Lambda \w^{22})$,  with $p(\w^{11}) p(\w^{22}) = p(\Lambda\w^{11}) p(\Lambda\w^{22})$. This will imply that
\begin{equation}\label{gammatoprove}
\sum_{\Gamma^1} p(\w^{11}) p(\w^{22})
=
\sum_{\Gamma^2} p(\w^{12}) p(\w^{21})
\end{equation}
from which $q = \det \mathbf{h}_{\bd A}(\mathbf{x},\mathbf{y})$ will follow.
To demonstrate the correspondence, the basic idea is to note that if $\w^{11} \in \curves_A(x^1,y^1)$ and $\w^{22} \in \curves_A(x^2,y^2)$ satisfy  $\LE(\w^{11}) \cap \w^{22} \neq \emptyset$, then there is a first site $v$ (sometimes called a pivot point) on $\LE(\w^{11})$ which is visited by the path $\w^{22}$. Now, consider the last time that $\w^{11}$ visits $v$ and the last time that $\w^{22}$ visits $v$. Interchanging the tails of the two excursions from the times of their respective last visits to $v$ produces two new excursions---one from $x^1$ to $y^2$ (written  $\Lambda\w^{11}$) and one from $x^2$ to $y^1$ (written $\Lambda\w^{22}$). Since the same sites are occupied in the new excursions as in the original excursions, we see that $p(\w^{11}) p(\w^{22}) = p(\Lambda\w^{11}) p(\Lambda\w^{22})$. Formally, we denote the two excursions 
$(\w^{11}, \w^{22}) \in \Gamma^1$ by $w^{11}:=[x^1, z_1, \ldots, z_m,y^1]$ and  $\w^{22}:=[x^2, w_1, \ldots, w_k,y^2]$, and we write $\LE(\w^{11})=[x^1,\hat{z_1}, \ldots, \hat{z_i},y^1]$.  Since $\LE(\w^{11}) \cap \w^{22} \neq \emptyset$,  we set $j_1$ to be the smallest positive integer $j$ such that $\hat{z}_j \in \{w_1, \ldots, w_k\}$ and $\ell_1$ to be the largest integer $\ell$ such that $z_{\ell}=\hat{z_{j_1}}$.  Finally, let $\ell_2$ be the largest integer $\ell$ such that $w_{\ell}=\hat{z_{j_1}}$. We now define
$\Lambda\w^{11}:=[x^1, z_1, \ldots, z_{\ell_1}, w_{\ell_2+1}, \ldots, w_{k}, y^2]$
and $\Lambda\w^{22}:=[x^2, w_1, \ldots, w_{\ell_2}, z_{\ell_1+1}, \ldots, z_{m}, y^1]$
so that $(\Lambda\w^{11}, \Lambda \w^{22}) \in \Gamma^2$ and $p(\w^{11}) p(\w^{22}) = p(\Lambda\w^{11}) p(\Lambda\w^{22})$. Conversely, if we consider $(\w^{12}, \w^{21}) \in \Gamma^2$, then necessarily $\LE(\w^{12}) \cap \w^{21} \neq \emptyset$. Interchanging the tails of the two excursions in the same manner just described therefore produces two new excursions---one from $x^1$ to $y^1$  (written $\Lambda^{-1}\w^{12}$) and one from $x^2$ to $y^2$ (written $\Lambda^{-1}\w^{21}$)---such that $\LE(\Lambda^{-1}\w^{12}) \cap \Lambda^{-1}\w^{21}\neq \emptyset$
and $p(\w^{12}) p(\w^{21}) = p(\Lambda^{-1}\w^{12}) p(\Lambda^{-1}\w^{21})$. This establishes the bijection between $\Gamma^1$ and $\Gamma^2$, establishes~(\ref{gammatoprove}) from which $q = \det \mathbf{h}_{\bd A}(\mathbf{x},\mathbf{y})$ follows, and completes the proof of the theorem.
\end{proof}

\begin{remark}
Fomin's identity (Theorem~\ref{fomintheorem}) as originally proved~\cite{Fom1} holds more generally for stationary Markov processes on discrete state spaces. The proof we present here of 
Theorem~\ref{FI_LERW_b} for the special case of two paths in $\Z^2$  holds with only minor modifications for two paths in $\Z^d$.
\end{remark}


\section{A brief introduction to SLE}\label{SLEintrosec}

The Schramm-Loewner evolution (SLE) is a one-parameter family of random growth processes introduced in 1999 by O.~Schramm~\cite{schramm} while considering possible scaling limits of loop-erased random walk. Since then a number of introductions to SLE have been written for a range of audiences. These include lecture notes by I.~Gruzberg~\cite{Gruzberg} and W.~Werner~\cite{WernerNotes}, and a book by G.~Lawler~\cite{SLEbook}. The purpose of this section (as the title suggests) is to provide a brief introduction. At times we will be a little casual sacrificing precision for intuition; the interested reader can find precise details in~\cite{SLEbook}.

Let $\Half = \{z \in \C : \Im(z) >0\}$ denote the upper half plane, and 
consider a simple (non-self-intersecting) curve $\gamma:[0,\infty) \to \overline{\Half}$ with $\gamma(0)=0$ and $\gamma(0,\infty) \subset \Half$.  For every fixed $t \ge 0$, the slit plane $\Halft := \Half \setminus \gamma(0,t]$ is simply connected and so by the Riemann mapping theorem, there exists a conformal transformation $g_t:\Halft \to \Half$.  The map $g_t$ is not unique, but we choose the unique one satisfying the hydrodynamic normalization $g_t(z) - z \to 0$ as $z \to \infty$.
It then follows that $g_t$ can be expanded as
\begin{equation}\label{expansion}
g_t(z) = z + \frac{b(t)}{z} + O \left(|z|^{-2}\right), \;\;\; z \to \infty,
\end{equation}
where $b(t) = \hcap(\gamma(0,t])$ is the \emph{half-plane capacity} of $\gamma$ up to time $t$. The half-plane capacity is related to how likely a Brownian motion starting from infinity is to hit the curve before hitting the real line $\R$. If $B_t$ is a two-dimensional Brownian motion, then
$$ \hcap(\gamma(0,t]) := \lim_{y \to \infty} y\, \E[\Im(B_{\tau}) | B_0= iy]$$
where $\tau$ is the first time that the Brownian motion hits either $\gamma(0,t]$ or $\R$.

\begin{figure}[h]
 \centering
 \includegraphics[height=1.3in]{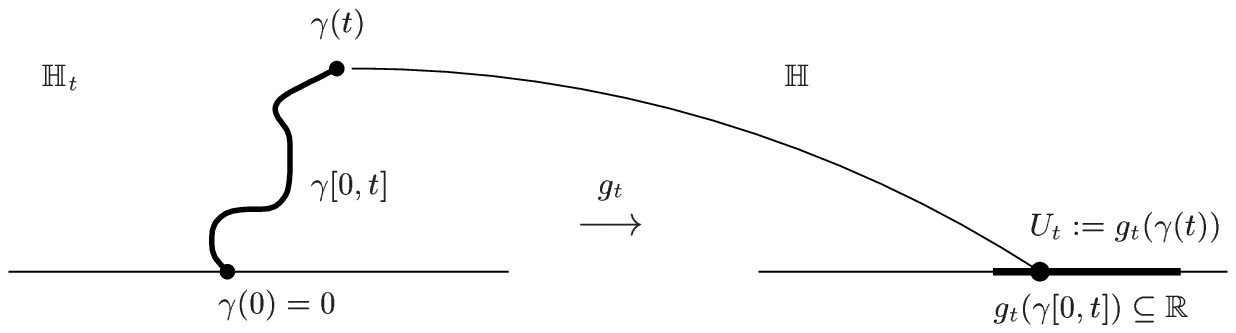}
\caption{The curve  $\gamma:[0,\infty) \to \overline{\Half}$ and the map $g_t:\Halft \to \Half$.}\label{SLE_map}
\end{figure}

For a slit plane such as $\Halft = \Half \setminus \gamma(0,t]$,  the map $g_t$ can be extended continuously to the boundary point $\gamma(t)$ of $\bd \Halft$.  With no additional assumptions on the simple curve $\gamma$, it can be shown that there is a unique point $U_t\in \R$ for all $t \ge 0$ with $U_t:= g_t(\gamma(t))$ and that the function $t \mapsto U_t$ is continuous. The notation is illustrated in Figure~\ref{SLE_map}.

The evolution of the curve $\gamma(t)$, or more precisely, the evolution of the conformal transformations $g_t: \Halft \to \Half$, can be described by a differential equation involving $U_t$. This is due to C.~Loewner who showed in 1923 that if
$\gamma$ is a curve as above such that its half-plane capacity $b(t)$ is $C^1$ and $b(t) \to \infty$ as $t \to \infty$, then for $z \in \Half$ with $z \not\in \gamma[0,\infty)$, the conformal transformations $\{g_t(z), t \ge 0\}$ satisfy the partial differential equation
\begin{equation}\label{LE}
\frac{\bd}{\bd t} \, g_t(z) = \frac{\dot b(t)}{g_t(z)-U_t}, \;\;\; g_0(z)=z.
\end{equation}
Note that if $b(t) \in C^1$ is an increasing function, then we can reparametrize the curve $\gamma$ so that $\hcap(\gamma(0,t]) = b(t)$.   This is the so-called parametrization by capacity and will be convenient for our purposes.

The obvious thing to do now is to start with a continuous function $t \mapsto U_t$  from $[0,\infty)$ to $\R$ and solve the Loewner equation~(\ref{LE}) for $g_t$.  Ideally, we would like to solve~(\ref{LE}) for  $g_t$, define simple curves
 $\gamma(t)$, $t \ge 0$, by setting $\gamma(t) = g_t^{-1}(U_t)$, and have $g_t$ map $\Half \setminus \gamma(0,t]$ conformally onto $\Half$. Although this is the correct intuition, it is not quite precise because we see from the denominator on the right-side of~(\ref{LE}) that problems can occur if $g_t(z)-U_t=0$.  Formally, if we let $T_z$ be the supremum of all $t$ such that the solution to~(\ref{LE}) is well-defined up to time $t$ with $g_t(z) \in \Half$, and we define $\Halft = \{ z : T_z>t\}$, then $g_t$ is the  unique conformal transformation of $\Halft$ onto $\Half$ with $g_t(z) - z \to 0$ as $t\to \infty$ and has expansion as in~(\ref{expansion}).

The novel idea of Schramm was to take the continuous function $U_t$ to be a one-dimensional Brownian motion starting at $0$ with variance parameter $\kappa \ge 0$. This leads to the following definition. The \emph{chordal Schramm-Loewner evolution with parameter $\kappa \ge 0$ with the standard parametrization} (or simply $\SLE_{\kappa}$)
is the random collection of conformal maps $\{g_t, \, t \ge 0\}$ obtained by solving the initial value problem
\begin{equation}\label{SLE}
\frac{\bd}{\bd t} \, g_t(z) = \frac{2}{g_t(z)-\sqrt{\kappa} W_t}, \;\;\; g_0(z)=z,
\end{equation}
where $W_t$ is a standard one-dimensional Brownian motion.

The question is now whether or not there exists a curve associated with the maps $g_t$. The answer is yes, although describing this curve requires the following deep theorem.

\begin{itemize}
\item If $0<\kappa \le 4$, then there exists a random simple curve $\gamma:[0,\infty) \to \overline{\Half}$ with $\gamma(0)=0$ and $\gamma(0,\infty) \subset \Half$. (That is, the curve  never re-visits $\R$.) Furthermore,  the maps $g_t$ obtained by solving~(\ref{LE}) are conformal transformations of  $\Half \setminus \gamma(0,t]$ onto $\Half$. For this range of $\kappa$, our intuition matches the theory!

\item For $4<\kappa < 8$, there exists a random curve $\gamma:[0,\infty) \to \overline{\Half}$. These curves have double points and they do hit $\R$, but they never cross themselves!  As such, $\Half \setminus \gamma(0,t]$ is not simply connected.  However, $\Half \setminus \gamma(0,t]$ does have a unique connected component containing $\infty$. This is $\Halft$ and the maps $g_t$ are conformal transformations of $\Halft$ onto $\Half$.  We think of $\Halft = \Half \setminus K_t$ where $K_t$ is the \emph{hull of $\gamma(0,t]$} visualized by taking $\gamma(0,t]$ and filling in the holes. In the case $0<\kappa \le 4$ where the curve is simple, we have $K_t = \gamma(0,t]$.

\item For $\kappa \ge 8$, there exists a random curve $\gamma:[0,\infty) \to \overline{\Half}$ which is space-filling! Furthermore, it has double points, but does not cross itself! As in the case $4 < \kappa <8$,  the maps $g_t$ are conformal transformations of $\Halft = \Half \setminus K_t$ onto $\Half$ where $K_t$ is the hull of $\gamma(0,t]$.
\end{itemize}

The case $\kappa \neq 8$ was established by S.~Rohde and O.~Schramm~\cite{RS} while the case $\kappa=8$ was proved by G.~Lawler, O.~Schramm, and W.~Werner~\cite{LSW}. As a result of this, we also refer to the curve $\gamma$ as chordal $\SLE_{\kappa}$. It is worth mentioning that SLE paths are extremely rough.  It has been shown by V.~Beffara~\cite{beffara} that the Hausdorff dimension of a chordal $\SLE_{\kappa}$ path is $\min\{1+\kappa/8, 2\}$. The Java applet simulation of SLE at
\begin{center}
\texttt{http://stat.math.uregina.ca/$\sim$kozdron/Simulations}
\end{center}
works particularly well for $0<\kappa<4$.

Since there exists a curve $\gamma$ associated with the maps $g_t$, it is possible to reparametrize it. As such, it can be shown that if $W_t$ is a standard one-dimensional Brownian motion, then  the solution to the initial value problem
\begin{equation}\label{SLEb}
\frac{\bd}{\bd t} \, g_t(z) = \frac{2/\kappa}{g_t(z)-W_t}, \;\;\; g_0(z)=z,
\end{equation}
is chordal $\SLE_{\kappa}$ parametrized so that  $\hcap(\gamma(0,t]) = 2t/\kappa$. For further  details on this point, see~\cite[Remark~6.7]{SLEbook}.

Finally, we would like to mention that chordal SLE as we have defined it can also be thought of as a probability measure on paths in the upper half plane $\Half$ connecting the boundary points $0$ and $\infty$. SLE is conformally invariant and so we can define chordal $\SLE_{\kappa}$ in any simply connected domain $D$ connecting distinct boundary points $z$ and $w$ to be the image of chordal $\SLE_{\kappa}$ in $\Half$ from $0$ to $\infty$ under a conformal transformation from $\Half$ onto $D$ sending $0 \mapsto z$ and $\infty \mapsto w$.

As already noted, O.~Schramm introduced SLE in 1999 while considering possible scaling limits of loop-erased random walk. Considerations from statistical mechanics suggested that the limit should be a random simple curve satisfying a type of conformal invariance property. Shortly thereafter, it was proved by G.~Lawler, O.~Schramm, and~W.~Werner~\cite{LSW} that, in fact, the scaling limit of loop-erased random walk can be given by SLE with parameter $\kappa=2$. (To give the precise technical details of the proof, they actually considered a slightly different version of SLE known as radial SLE which is a measure on paths connecting an interior point to a boundary point. The extension of showing convergence of loop-erased random walk to chordal $\SLE_2$ will appear as part of the Ph.D.\ dissertation of F.~Johansson of KTH Stockholm.)


\section{Review of the excursion Poisson kernel}\label{epksect}
 
The excursion Poisson kernel is formally defined as the normal derivative of the usual Poisson kernel. However, it is also the mass of the Brownian excursion measure (which itself is the scaling limit of simple random walk excursion measure), and the original motivation for studying the excursion Poisson kernel was in this context. Further details may be found in~\cite{Koz} and~\cite{KozL}.

 Suppose that $D \subset \C$ is a simply connected Jordan domain and that $\bd D$ is locally analytic at $x$ and $y$. The \emph{excursion Poisson kernel} is defined as
 \[
 H_{\bd D}(x,y) := \lim_{\eps \to 0} \frac{1}{\eps} \, H_D(x + \eps \n_x, y)
 \]
 where $H_D(z,y)$ for $z \in D$ is the usual Poisson kernel, and $\n_x$ is the unit normal at $x$ pointing into $D$. The excursion Poisson kernel satisfies the following important conformal covariance property; see~\cite[Proposition~2.11]{KozL}. If $f:D \to D'$ is a conformal transformation where $D' \subset \C$ is also a simply connected Jordan domain, and $\bd D'$ is locally analytic at $f(x)$, $f(y)$, then
\begin{equation}\label{EPKconfcov}
H_{\bd D}(x,y) = |f'(x)| |f'(y)| H_{\bd D'}(f(x),f(y)).
\end{equation}
Explicit formul\ae\ are known when $D=\D$, the unit disk, or $D=\Half$, namely
 \begin{equation*}
  H_{\bd \D}(x,y) = \frac{1}{\pi \, | y-x|^2} =\frac{1}{2\pi(1-\cos(\arg y - \arg x))} \end{equation*}
and
\begin{equation}\label{EPKform}
H_{\bd \Half}(x,y) = \frac{1}{\pi(y-x)^2}.
\end{equation}
Suppose now that $x^1, \ldots, x^n, y^1, \ldots, y^n$ are distinct boundary points at which $\bd D$ is locally analytic, 
let $f: D \to D'$ be a conformal transformation, 
and assume that $\bd D'$ is also locally analytic at  $f(x^1), \ldots, f(x^n), f(y^1), \ldots, f(y^n)$.
It follows~\cite[Proposition~2.16]{KozL} that if $\mathbf {H}_{\bd D}(\x,\y) := [H_{\bd D}(x^i, y^{\ell})]_{1 \le i, \ell \le n}$ denotes the $n \times n$ \emph{hitting matrix}
\[
\mathbf {H}_{\bd D}(\x,\y) := 
\begin{bmatrix}
H_{\bd D}(x^1, y^1)   &\cdots    & H_{\bd D}(x^1, y^n) \\
\vdots                         &\ddots    &\vdots \\
H_{\bd D}(x^n, y^1)   &\cdots    &H_{\bd D}(x^n, y^n) \\
\end{bmatrix}
 \]
 then
 \begin{equation}\label{feb20eq2}
 \det  \mathbf {H}_{\bd D}(\x,\y) = \left( \prod_{j=1}^n |f'(x^j)| \; |f'(y^j)| \right) \det  [H_{\bd D'}(f(x^i), f(y^\ell))]_{1 \le i, \ell \le n}. 
 \end{equation}
It now follows from~(\ref{EPKconfcov}) and~(\ref{feb20eq2}) that
\begin{equation}\label{EPKratio}
\frac{
\det  \mathbf {H}_{\bd D}(\x,\y)}{
\ds \prod_{i=1}^n {H}_{\bd D}(x^i,y^i)}
\end{equation}
is a conformal \emph{invariant}.

 It is worth noting that $H_{\bd D}(x,y)$ can be defined even if $\bd D$ is \emph{not} locally analytic at $x$, $y$. Simply let $f:\D \to D$ and take
$$H_{\bd D}(x,y) := |f'(f^{-1}(x))|^{-1}|f'(f^{-1}(y))|^{-1} H_{\bd \D}(f^{-1}(x),f^{-1}(y)).$$

The reader will no doubt notice the similarities between the Brownian motion functional~(\ref{EPKratio}) and the simple random walk functional~(\ref{conditionalFomin}). The first approach to establishing a scaling limit of Fomin's identity~\cite{KozL} involved showing that an appropriate limit of~(\ref{conditionalFomin}) existed as the lattice spacing $\delta \to 0$. In fact, as an extension of that work, it is shown in~\cite{Koz} that simple random walk excursion measure converges to Brownian excursion measure on any simply connected domain with Jordan boundary.


\section{The non-intersection probability of $\SLE_2$ and Brownian motion}\label{nonintsec}

Since Fomin's identity allows us to calculate the probability that loop-erased random walk and simple random walk do not intersect, the natural continuous analogue is the probability that $\SLE_2$  and Brownian motion  do not intersect. 

\begin{figure}[h]
 \centering
 \includegraphics[height=1.5in]{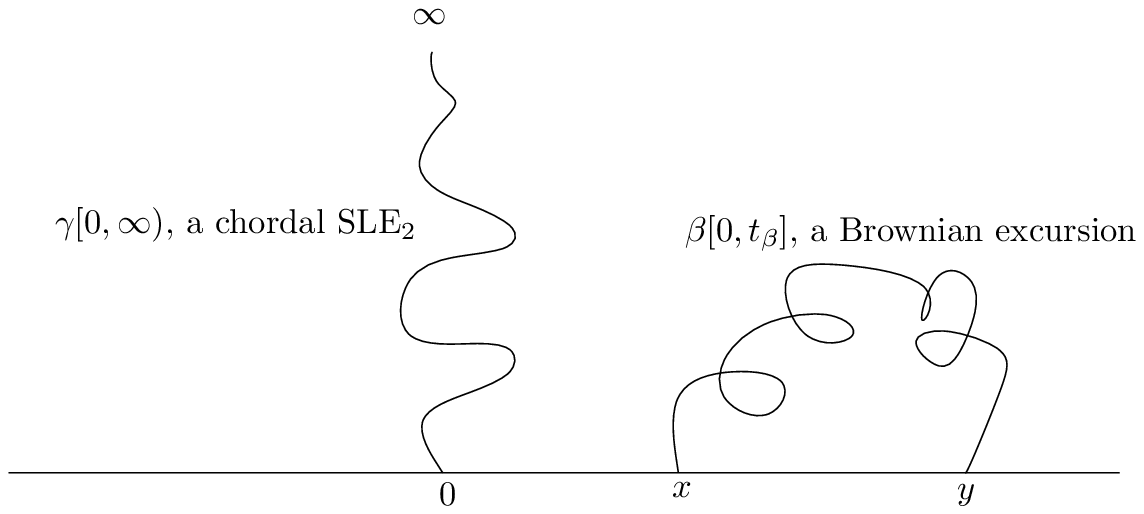}
\caption{Schematic representation of $\Prob\{ \gamma[0,\infty) \cap \beta[0,t_{\beta}] = \emptyset \}$. }\label{SLEfig}
\end{figure}

Suppose that $\gamma:[0,\infty) \to \overline{\Half}$ is a chordal $\SLE_2$ from $0$ to $\infty$ in $\Half$. Suppose further that $0<x<y<\infty$ are real numbers and let $\beta:[0,t_{\beta}] \to \overline{\Half}$ be a Brownian excursion from $x$ to $y$ in $\Half$. Hence, our goal is to determine
$\Prob\{ \gamma[0,\infty) \cap \beta[0,t_{\beta}] = \emptyset \}$
and show that it can be expressed in terms of the determinant of the Brownian excursion hitting matrix; see Figure~\ref{SLEfig}.  Notice the similarity between the following theorem which may be called Fomin's identity for $\SLE_2$ and Theorem~\ref{FI_LERW_b}. The proof we include is an expanded version (giving all the details) of the one in~\cite{KozL2}.

 \begin{thrm}[Fomin's Identity for $\SLE_2$]\label{FI_SLE}
If $x$, $y \in \R$ with $0<x<y<\infty$ and that $\gamma$, $\beta$ are as above, then 
\begin{equation}\label{scalefomin}
 \Prob\{ \,\gamma[0,\infty) \cap \beta[0,t_{\beta}] = \emptyset\, \}
=
\frac{\det  \mathbf {H}_{\bd \D}(f(\x),f(\y)) }{ 
H_{\bd \D}(f(0), f(\infty)) \, H_{\bd \D}(f(x), f(y))
}
\end{equation}
where $f: \Half \to \D$ is a conformal transformation.
\end{thrm}
 
Our strategy for establishing this result will be as follows.  We will first determine an explicit expression for $\Prob\{ \,\gamma[0,\infty) \cap \beta[0,t_{\beta}] = \emptyset\, \}$,  and we will then show that this explicit expression is the same as the right side of~(\ref{scalefomin}). 

\begin{proof}
For every $0 < t < \infty$, let $\Halft$ denote the slit-plane
$\Halft = \Half \setminus \gamma(0,t]$ which implies that 
$$\Prob\{ \gamma[0,t] \cap \beta[0,t_{\beta}] = \emptyset  \} = \E \left[ \frac{H_{\bd \Halft}(x,y)}{H_{\bd \Half}(x,y)}\right].$$
Therefore, letting $t \to \infty$ we conclude that
\begin{equation}\label{eq3}
\Prob\{ \gamma[0,\infty) \cap \beta[0,t_{\beta}] = \emptyset  \} =  \lim_{t \to \infty}\E \left[ \frac{H_{\bd \Halft}(x,y)}{H_{\bd \Half}(x,y)}\right]
=  \E \left[ \lim_{t \to \infty} \frac{H_{\bd \Halft}(x,y)}{H_{\bd \Half}(x,y)}\right].
\end{equation}
Let $g_t:\Halft \to \Half$ be the unique conformal transformation satisfying the hydrodynamic normalization $g_t(z)-z = o(1)$ as $z \to \infty$. As indicated in Section~\ref{SLEintrosec}, it is well-known that $g_t$ satisfies the chordal Loewner equation, namely
\begin{equation}\label{SLEeq}
\frac{\bd}{\bd t} \, g_t(z) = \frac{1}{g_t(z) +W_t}, \;\;\; g_0(z)=z,
\end{equation}
where $W_t$ is a standard Brownian motion. (This follows from~(\ref{SLEb}) by noting that if $W_t$ is a standard Brownian motion, then so too is $-W_t$.)
We now map $\Halft$ to $\Half$ by $g_t$ and use
conformal covariance~(\ref{EPKconfcov}) to conclude that
$H_{\bd \Halft}(x,y) = g_t'(x)g_t'(y) H_{\bd \Half}(g_t(x),g_t(y))$
and so
\begin{equation}\label{eq1feb20}
\frac{H_{\bd \Halft}(x,y)}{H_{\bd \Half}(x,y)}=  \frac{g_t'(x)g_t'(y) H_{\bd \Half}(g_t(x),g_t(y))}{H_{\bd \Half}(x,y)} =  (y-x)^2 \cdot \frac{g_t'(x)g_t'(y)}{(g_t(y)-g_t(x))^2}
\end{equation}
where the last equality follows from the explicit form of $H_{\bd \Half}$ in~(\ref{EPKform}).
Let 
$$J_t :=  \frac{g_t'(x)g_t'(y)}{(g_t(y)-g_t(x))^2} \;\;\; \text{and set} \;\;\; 
J_{\infty} := \lim_{t\to \infty}J_t.$$
To be consistent with notation in other papers (such as~\cite{KozL2}), let
$\tilde H^{*}(x,y) := \Prob\{ \gamma[0,\infty) \cap \beta[0,t_{\beta}] = \emptyset  \}$
so that~(\ref{eq3}) and~(\ref{eq1feb20}) give
\begin{align}\label{eq2}
\tilde H^{*}(x,y) 
= (y-x)^2\E \left[ \lim_{t \to \infty} \frac{g_t'(x)g_t'(y)}{(g_t(y)-g_t(x))^2} \right]
&= (y-x)^2\E \left[ \lim_{t \to \infty} J_t\right] \notag\\
&= (y-x)^2\E [J_{\infty}].
\end{align}
Since our goal is to compute
$\tilde H^{*}(x,y)$
we will derive a differential equation for $\tilde H^{*}(x,y)$.  Let $X_t := g_t(x)+W_t$ and $Y_t := g_t(y) + W_t$ where $g_t$ and $W_t$ are as in~(\ref{SLEeq}) so that 
$$dX_t = \frac{1}{X_t} dt + dW_t \;\;\; \text{and} \;\;\; dY_t = \frac{1}{Y_t} dt + dW_t.$$
Some routine calculations give
\[     \frac{\bd}{\bd t} \, [\log g_t'(x)] = -\frac{1}{X_t^2}, \;
              \frac{\bd}{\bd t} \,[\log g_t'(y)] = - \frac{1}{Y_t^2}, \; \text{and} \;
      \frac{\bd}{\bd t} \,[\log(g_t(y) - g_t(x))]  =
   -\frac{1}{X_t \, Y_t},  \]
and so we see that
\begin{align*}
J_t = \exp\left\{ \log J_t \right\} 
&=J_0 \, \exp\left\{\int_0^t \frac{\bd}{\bd s} \; [\log J_s] \; ds \right\}\\
&=\frac{1}{(y-x)^{2}} \, \exp\left\{- \int_0^t \left( \frac{1}{X_s} - \frac{1}{Y_s}
\right)^2 \; ds \right\}
\end{align*}
since
$$J_0 = \frac{g_0'(x)g_0'(y)}{(g_0(y)-g_0(x))^2} =\frac{1}{(y-x)^2}.$$
Hence~(\ref{eq2}) implies that
\[    \tilde H^*(x,y) = \E \left[
 \exp\left\{- \int_0^{\infty} \left( \frac{1}{X_s} - \frac
   1 {Y_s} \right)^2 \; ds \right\}\right].\]
It now follows from the (usual) Markov property that $J_t \, \tilde H^*(X_t,Y_t)$ is a 
martingale. That is, 
if $M_t := \E[J_{\infty}| \mathcal{F}_t]$ so that $M_t$ is a martingale, then
\begin{align*}
&M_t 
=\E \left[ \frac{1}{(y-x)^2} \;\exp\left\{- \int_0^{\infty} \left( \frac{1}{X_s} - \frac{1}{Y_s}
 \right)^2 \; ds \right\} 
\; \bigg| \;\mathcal{F}_t\;\right]\\
&=\frac{1}{(y-x)^2}\exp\left\{- \!\int_0^{t} \left( \frac{1}{X_s} - \frac{1}{Y_s} \right)^2  ds \right\} 
\E \left[  \exp\left\{- \!\int_{t}^{\infty} \left( \frac{1}{X_s} - \frac{1}{Y_s}
 \right)^2  ds \right\} 
 \bigg| \mathcal{F}_t\right]\\
&=J_t\, \tilde H^*(X_t,Y_t).
\end{align*}
It\^o's formula at $t=0$ now implies that 
\begin{equation}\label{ItoEqn}
  -\left(\frac 1x - \frac 1y \right)^2 \, \tilde H^*  + \frac 1x \, 
    \frac{\bd \tilde H^*}{\bd x}   
+ \frac 1y \, \frac{\bd \tilde H^*}{\bd y} + \frac 12 \, \frac{\bd^2 \tilde H^*}{\bd x^2} 
     + \frac 12 \, \frac{\bd^2 \tilde H^*}{\bd y^2}
   + \frac{\bd^2 \tilde H^*}{\bd x\bd y}=0
\end{equation}
Since the probability in question only depends on the ratio $x/y$, we see that
$\tilde H^*(x,y) = \phi(x/y)$ for some function $\phi$. Thus, we find
\[   \frac{\bd \tilde H^*}{\bd x} = y^{-1} \, \phi'(x/y), \;\;\; \frac{\bd \tilde H^*}{\bd y} = - x \, y^{-2} \, \phi'(x/y), \;\;\;
   \frac{\bd^2 \tilde H^*}{\bd x^2}= y^{-2} \, \phi''(x/y), \]
\[   \frac{\bd^2 \tilde H^*}{\bd y^2} = 2\,  x \, y^{-3} \, \phi'(x/y) + x^2 \, y^{-4} \,   \phi''(x/y), \;\;\;
 \frac{\bd^2 \tilde H^*}{\bd x\bd y}= - y^{-2} \, \phi'(x/y) - x \, y^{-3} \, \phi''(x/y),\]
so that after substituting into~(\ref{ItoEqn}), multiplying by $y^2$, letting $u=x/y$,
and combining terms, we have
\begin{equation}\label{hypergeomFomin}
 u^2 \, (1-u)^2 \,\phi''(u)  +  2 \,
 u \, (1-u) \, \phi'(u) -2 (1-u)^2 \, \phi(u) =0.
\end{equation}
Observe, however, that~(\ref{hypergeomFomin}) is equivalent to
\begin{equation}\label{hypergeomFomin2}
 u^2 \, (1-u) \,\phi''(u)  +  2 \,
 u \, \phi'(u) -2 (1-u) \, \phi(u) =0
\end{equation}
since $0<u<1$. 
The second-order ordinary differential equation~(\ref{hypergeomFomin2}) has regular singular points at $0$, $1$, and $\infty$, and so we know that it is possible to 
transform it into a hypergeometric differential equation.
By writing~(\ref{hypergeomFomin2}) as 
\begin{equation}\label{hypergeomFomin2-pform}
\phi''(u)  + \left[ \frac{2}{u} - \frac{2}{u-1} \right] \phi'(u) + \left[ \frac{2}{u^2(u-1)} - \frac{2}{u(u-1)} \right]  \phi(u) =0
\end{equation}
we see that we have a case of Riemann's differential equation whose complete set of solutions (see~(15.6.1) and~(15.6.3) of~\cite{AbSteg}) can be denoted by Riemann's $P$-function
$$\phi(u) = P 
\left\{
\begin{matrix}
0 &\infty &1 &\\
1 &-2 &3 &u\\
-2 &1 &0 &\\
\end{matrix}
\right\}.$$
By now considering~(15.6.11) of~\cite{AbSteg}, the transformation formula for Riemann's $P$-function for reduction to the hypergeometric function, we see that the appropriate change-of-variables to apply is  $\psi(u) := u^{-1} (1-u)^{-3} \phi(u)$ noting that this is permitted by the constraint $0<u<1$.
 Thus,~(\ref{hypergeomFomin2}) implies
\begin{equation}\label{hypergeomFomin3}
 u \, (1-u) \,\psi''(u)  +  (4-8u) \, \psi'(u) -10\, \psi(u) =0.
\end{equation}
We see that~(\ref{hypergeomFomin3}) is now a well-known hypergeometric differential equation~\cite{AbSteg} whose general solution is given by
$$\psi(u) =C_1 \frac{2-u}{(1-u)^{3}} + C_2 \frac{1-2u}{u^3(1-u)^3}.$$
This implies that the general solution to~(\ref{hypergeomFomin2})
is
$\phi(u)=
C_1 u(2-u) + C_2 u^{-2}(1-2u)$.
However, physical considerations dictate that $\phi(u) \to 0$ as $u \to 0+$ and $\phi(u) \to 1$ as $u \to 1-$, and so $C_2=0$ and $C_1=1$.
Thus, $\phi(u)= u(2-u)$ and so we find
\begin{equation}\label{fominprob}
\Prob\{ \gamma[0,\infty) \cap \beta[0,t_{\beta}] = \emptyset  \} =
\tilde H^*(x,y) = \phi(x/y) =\frac{x}{y}\left(2-\frac{x}{y}\right).
\end{equation}
As already noted, the probability in question only depends on the ratio $x/y$, and so it suffices without loss of generality to assume that $0<x<1$ and $y=1$. Furthermore,
we may assume that the conformal transformation $f: \Half \to \D$ is given by
\begin{equation}\label{conformmap}
f(z) = \frac{iz+1}{z+i},
\end{equation}
so that $f(0)=-i$, $f(y)=f(1)=1$, $f(\infty)=i$, and
\[
f(x) = \left(\frac{2x}{x^2+1} \right) + i \left( \frac{x^2-1}{x^2+1} \right) = \exp \left\{ -i \arctan \left( \frac{1-x^2}{2x} \right) \right\}.
\]
Writing $f(x) = e^{i \theta}$, we find that
\begin{align*}
&\frac{\det  \mathbf {H}_{\bd \D}(f(\x),f(\y)) }{ 
H_{\bd \D}(f(0), f(\infty))  H_{\bd \D}(f(x), f(y)) }
= \frac{
H_{\bd \D}(-i,i)  H_{\bd \D}(e^{i\theta},1)-H_{\bd \D}(-i,1)  H_{\bd \D}(e^{i\theta},i)}
{H_{\bd \D}(-i,i)  H_{\bd \D}(e^{i\theta},1) }\\
&\qquad\qquad\qquad\qquad= \frac{ \frac{1}{2\pi (1-\cos \pi)} \frac{1}{2\pi (1-\cos \theta)}
-\frac{1}{2\pi (1-\cos (\frac{\pi}{2}) )} \frac{1}{2\pi (1-\cos(\frac{\pi}{2}+\theta))}}
{\frac{1}{2\pi (1-\cos \pi)} \frac{1}{2\pi (1-\cos \theta)}}\\
&\qquad\qquad\qquad\qquad= \frac{2\cos \theta + \sin \theta - 1}{1 + \sin \theta}.
\end{align*}
Since 
$\theta = - \arctan \left(\frac{1-x^2}{2x} \right)$
we see that
$\cos \theta  = \frac{2x}{x^2+1}$  and $\sin \theta = \frac{1-x^2}{x^2+1}$
which upon substitution gives
\[ 
\frac{2\cos \theta + \sin \theta - 1}{1 + \sin \theta}
= \frac{ \frac{4x}{x^2+1} + \frac{1-x^2}{x^2+1} - 1}{1 +  \frac{1-x^2}{x^2+1}}
= \frac{4x-2x^2}{2} = x(2-x).\]
Comparison with~(\ref{fominprob}) now yields the result, and the theorem is proved.
\end{proof}      

\begin{example}
Suppose that $\gamma$ is a chordal $\SLE_2$ from $0$ to $\infty$ in $\Half$, and let $\beta$ be a Brownian excursion from $x=1/2$ to $y=1$ in $\Half$. Let 
$f(z)$ be as in~(\ref{conformmap})
which is a conformal transformation of $\Half$ onto $\D$ with $f(0) = i$, $f(1)=1$, $f(\infty)=-i$. Also notice that $f(1/2) =\exp\{-i \arctan(3/4)\}$. A simple calculation then shows that
$$\Prob\{ \gamma[0,\infty) \cap \beta[0,t_{\beta}] = \emptyset  \}= \frac{2\cdot\frac{4}{5} + \frac{3}{5} - 1}{1 + \frac{3}{5}} = \frac{1}{2}\left(2-\frac{1}{2}\right) = \frac{3}{4}.$$
\end{example}

\begin{remark}
As the reader has no doubt discovered, by working in $\Half$ and $\D$ it is possible to perform explicit calculations. Since the quantity on the right side of~(\ref{scalefomin}) is known to be a conformal invariant as in~(\ref{EPKratio}), we can show, with a combination of conformal transformations, that the probability a chordal $\SLE_2$ avoids a Brownian excursion in any simply connected domain $D$ is given by the appropriate determinant of the matrix of excursion Poisson kernels.
\end{remark}

\begin{cor}
Suppose that $D \subset \C$ is a bounded, simply connected planar domain, and that $x^1, x^2, y^2, y^1$ are 
four points ordered counterclockwise around $\bd D$. The probability a chordal $\SLE_2$ from $x^1$ to $y^1$ in $D$ does not intersect a Brownian excursion from $x^2$ to $y^2$ in $D$ is 
$ \Phi(x^2) \left(2- \Phi(x^2) \right)$
where $\Phi: D \to \Half$ is the conformal transformation with $\Phi(x^1) =0$, $\Phi(y^1)=\infty$, $\Phi(y^2)=1$.
\end{cor}

This statement can be easily modified to cover the case when $D$ is unbounded and/or the case when $\infty$ is one of the boundary points.

\begin{remark}
In the case of $n\ge 2$ paths, it is shown in~\cite{KozL} that the scaling limit of the determinant of the simple random walk hitting matrix~(\ref{conditionalFomin}) is the  determinant of the Brownian excursion hitting matrix~(\ref{EPKratio}). The proof of this result does not employ any SLE techniques. A formula analogous to Theorem~\ref{FI_SLE} relating the determinant of the Brownian excursion hitting matrix to a probability involving $n>2$ chordal $\SLE_2$ paths and Brownian excursions has not yet appeared, although there has been work done constructing a finite measure on $n \ge 2$ mutually avoiding $\SLE_2$ paths (see~\cite{Dub} and~\cite{KozL}) in which  the determinant of the Brownian excursion hitting matrix is related to the mass of this configurational measure.
\end{remark}

\proof[Acknowledgements]
The author would like to express his gratitute to the Pacific Institute for the Mathematical Sciences at the University of British Columbia in Vancouver, BC. PIMS hosted a February 2007 visit by the author during which the preliminary draft of this paper was completed. The author also wishes to thank Prof.~G.~Lawler in collaboration with whom much of this work was originally done; see~\cite{KozL} and~\cite{KozL2}. Special thanks are owed to the anonymous referee who provided several extremely useful suggestions for improving the exposition.


\end{document}